\documentclass[12pt,reqno]{amsart}
\usepackage{amsfonts}
\usepackage{amssymb,color}
\usepackage{lipsum}
\usepackage{mathrsfs}
\usepackage[colorlinks=true]{hyperref}
\hypersetup{urlcolor=blue, citecolor=blue}

\usepackage[margin=3cm, a4paper]{geometry}

\newtheorem{theo}{Theorem}[section]

\newtheorem{defi}[theo]{Definition}
\newtheorem{coro}[theo]{Corollary}

\newtheorem{rema}[theo]{Remark}
\numberwithin{equation}{section}

\newcommand{\bal}{\begin{align}}
\newcommand{\bbal}{\begin{align*}}
\newcommand{\beq}{\begin{equation}}
\newcommand{\eeq}{\end{equation}}
\newcommand{\bca}{\begin{cases}}
\newcommand{\eca}{\end{cases}}

\newcommand{\pa}{\partial}
\newcommand{\fr}{\frac}
\newcommand{\na}{\nabla}
\newcommand{\De}{\Delta}

\newcommand{\cd}{\cdot}
\newcommand{\ep}{\varepsilon}
\newcommand{\dd}{\mathrm{d}}

\newcommand{\R}{\mathbb{R}}

\newcommand{\D}{\mathrm{div}}

\begin{document}

\subjclass[2010]{}
\footnotetext{*Corresponding author.}

\title[Navier-Stokes equations with the Coriolis force]{Global large solutions for the Navier-Stokes equations with the Coriolis force}

\author[J. Li]{Jinlu Li}
\address{School of Mathematics and Computer Sciences, Gannan Normal University, Ganzhou 341000, P. R. China}
\email{lijinlu@gnnu.cn}

\author[J. Sun]{Jinyi Sun*}
\address{College of Mathematics and Statistics, Northwest Normal University, Lanzhou 730070, P. R. China}
\email{sunjy@nwnu.edu.cn}

\author[M. Yang]{Minghua Yang}
\address{Department of Mathematics, Jiangxi University of Finance and Economics, Nanchang 330032, P. R. China}
\email{ymh20062007@163.com}

\begin{abstract}
In this paper, we construct a class of global large solution to the three-dimensional Navier-Stokes equations with the Coriolis force in critical Fourier-Besov space $\dot{FB}^{2-\frac{3}{p}}_{p,r}(\mathbb{R}^3)$. In fact, our choice of special initial data $u_0$ can be arbitrarily large in $\dot{FB}^{s}_{p,r}(\mathbb{R}^3)$ for any $s\in\R$ and $1\leq p,r\leq \infty$.
\end{abstract}

\maketitle

\section{Introduction and main result}

Rotating fluid equations have important applications in meteorology and oceanography, particularly in the models describing large-scale ocean and atmosphere
flows. The Coriolis force, arising from the rotation of the Earth, plays a significant
role in such systems.

In 1868, Kelvin first observed that a sphere, moving along the axis of uniformly rotating water, takes with it a column of liquid as if this were a rigid mass, and pioneered the research on the motion of rotating fluid, see \cite{Fuentes}. Later on, Taylor\cite{Taylor} and Proudman\cite{Proudman} strictly proved that high-speed rotation brings about a vertical rigidity in the fluid described by the Taylor-Proudman theorem: Under a fast rotation the velocity of all particles located on the same vertical line is horizontal and constant.

%From the mathematical viewpoint,

Mathematically, the Coriolis forces give rise to the so-called
Poincar\'{e} waves, which are dispersive waves. Poincar\'{e} waves propagate in both directions with extremely fast speed in the propagation domain, and the waves with different wavenumbers move at different speeds. This makes that the nonlinear interactions between different modes are typically less significant.

On the other hand, Poincar\'{e} waves are a kind of high frequency wave, whose particles not only have vibrations parallel to the propagation direction, but also have vibrations perpendicular to the propagation direction. Therefore, one of the major difficulties encountered in understanding dynamics of rotating fluid is the influence of the oscillations generated by Coriolis forces.

In the paper, we consider the following Cauchy problem of three-dimensional incompressible Navier-Stokes equations with the Coriolis force:
\bal\label{NSC}\begin{cases}
\partial_tu-\De u+\Omega e_3\times u+u\cd\na u+\na p=0,\\
\D u=0, \quad u(0,x)=u_0,
\end{cases}\end{align}
where the unknown functions $u=(u_1, u_2, u_3)$ and $p$ denote velocity field and pressure, respectively; $\Omega\in \mathbb{R}$ is the Coriolis parameter, which is twice angular velocity of the rotation around the vertical unit vector $e_3=(0,0,1)$, and $\Omega e_3\times u$ represents the so-called Coriolis force; $u_0$ is the given initial velocity.

The behavior of fluid flows in rapidly rotating environments is fundamentally different from that of non-rotating flows.

When $\Omega=0$, \eqref{NSC} reduces to the problem of classical three-dimensional incompressible Navier-Stokes equations, which have been widely studied during the past seventy
years. It has been proved that \eqref{NSC} with $\Omega=0$ is globally well-posed for small initial data, see \cite{Cannone 1994,Fujita,Kato,Koch,Kukavica,Lei 2011}. For more results of large initial data  with special structures in various scaling invariant spaces which generate unique global solutions to \eqref{NSC}, we refer the reader to see \cite{Chemin2006-2,Chemin2010,Chemin2013,Chemin2015,Lei 2015,Li2019,Liu2018} and the references therein. We note that the global regularity or global well-posedness issue of the three-dimensional incompressible Navier-Stokes equations for arbitrarily large initial data is still a challenging open problem.

When $\Omega\neq0$, it is a very remarkable fact that
\eqref{NSC} admits a global solution for arbitrary large initial data, provided that the speed of rotation is fast enough. More precisely, when $\Omega$ is large enough, by taking the full exploration of dispersive effects of the Coriolis forces, the existence and uniqueness of global solution has been proved for the periodic large data in \cite{Babin 1997,Babin 1999}, for the spatially almost periodic large data in \cite{Yoneda 2011}, and for the decay large data see \cite{Chemin 2006, Iwabuchi 2013, Koh, Sun 2017}. For any given and fixed $\Omega$, we refer to \cite{Giga, Hieber, Iwabuchi 2014, Konieczny} for the global well-posedness of \eqref{NSC} with uniformly small initial data $u_0$. Especially, it has been proved that \eqref{NSC} is globally well-posed for small initial data in $\dot{FB}^{2-\frac{3}{p}}_{p,r}(\mathbb{R}^3)$ ($1< p,r\leq\infty$) and $\dot{FB}^{-1}_{1,r}(\mathbb{R}^3)$ ($1\leq r\leq 2$), and is ill-posed in  $\dot{FB}^{-1}_{1,r}(\mathbb{R}^3)$ ($r>2$), see \cite{Iwabuchi 2014,Konieczny}.

It is now a natural question to ask whether there exists a unique global solution to \eqref{NSC} for any given and fixed $\Omega$, if the initial data is not small in $\dot{FB}^{2-\frac{3}{p}}_{p,r}(\mathbb{R}^3)$ ($1\leq p,r\leq\infty$). Based on a full understanding of the structure of the equation \eqref{NSC}, we shall prove that \eqref{NSC} is globally well-posed for some special initial data $u_0$ whose $\dot{FB}^{2-\frac{3}{p}}_{p,r}(\mathbb{R}^3)$-norm can be arbitrarily large, namely, $||u_0||_{\dot{FB}^{2-\frac{3}{p}}_{p,r}}\gg 1$, for any $1\leq p,r\leq \infty$.

We first recall the definition of the Fourier-Besov spaces $\dot{FB}^s_{p,r}(\mathbb{R}^3)$.
As usual we denote by $\mathscr{S}(\mathbb{R}^3)$
the space of Schwartz functions on $\mathbb{R}^3$, and by $\mathscr{S}'(\mathbb{R}^3)$ the space of tempered distributions on $\mathbb{R}^3$.
Choose radial function $\psi\in\mathscr{S}(\mathbb{R}^3)$ such that its Fourier transform $\hat{\psi}$ satisfies the following properties:
$$\textrm{supp}~\hat{\psi}\subset \mathcal{C}:=\{\xi\in \mathbb{R}^3:\frac{3}{4}\leq|\xi|\leq\frac{8}{3}\},
$$
and
$$\sum_{j\in\mathbb{Z}} \hat{\psi}(2^{-j}\xi)=1 ~\quad \textrm{for all } \xi\in \mathbb{R}^3\setminus\{0\}.$$
Let $\psi_j(x):=2^{3j}\psi(2^jx)$ for $j\in\mathbb{Z}$ and $\mathscr{S}'_h(\mathbb{R}^3):=\mathscr{S}'(\mathbb{R}^3)/\mathcal{P}[\mathbb{R}^3]$, where $\mathcal{P}[\mathbb{R}^3]$ denotes the linear space of polynomials on $\mathbb{R}^3$. The homogeneous dyadic blocks $\Delta_j$ is defined by
$$\Delta_j f:=\psi_j \ast f$$  $\textrm{for}~j\in \mathbb{Z}\ \ \text{and}\ \ f\in \mathscr{S}'(\mathbb{R}^3)$.
Then the Fourier-Besov spaces $\dot{FB}^s_{p,r}(\mathbb{R}^3)$ are defined as follows:
\begin{defi}
For $s\in\mathbb{R}$ and $1\leq p,r\leq \infty$, the Fourier-Besov space $\dot{FB}^s_{p,r}(\mathbb{R}^3)$ is defined to be the
set of all tempered distributions $u\in \mathscr{S}'_h(\mathbb{R}^3) $ such that
$$\|u\|_{\dot{FB}^{s}_{p,r}}:=\Big\|\Big\{2^{js}\|
\widehat{\Delta_ju}\|_{L^p}\Big\}_{j\in\mathbb{Z}}\Big\|_{\ell^r(\mathbb{Z})}<\infty.$$
\end{defi}

\begin{rema}\label{re1.2}
It is easy to show that $||u||_{\dot{FB}^{0}_{p,1}}=||\hat{u}||_{L^p}$.
\end{rema}

Let $U$ satify the following linear system:
\bbal\begin{cases}
\pa_tU-\De U+\Omega e_3\times U+\na p'=0,
\\\D U=0, \quad U(0,x)=u_0.
\end{cases}\end{align*}
According to \cite{Hieber}, we can show that $U$ have the following explicit form:
\begin{equation}\label{1.2}
\hat{U}=\cos(\Omega\frac{\xi_3}{|\xi|}t)e^{-|\xi|^2t}\hat{u}_0
+\sin(\Omega\frac{\xi_3}{|\xi|}t)\frac{1}{|\xi|}e^{-|\xi|^2t}
\begin{pmatrix}
\xi_3\hat{u}_{0}^2-\xi_2 \hat{u}_{0}^3\\ -\xi_3\hat{u}_{0}^1+\xi_1 \hat{u}_{0}^3\\ \xi_2\hat{u}_{0}^1-\xi_1\hat{u}_{0}^2
\end{pmatrix},
\end{equation}
and it is easy to check that
$$||U||_{L^\infty(0,\infty;\dot{FB}^{-1}_{1,1})}+
||U||_{L^1(0,\infty;\dot{FB}^{1}_{1,1})}\leq C||u_0||_{\dot{FB}^{-1}_{1,1}}.$$
The main result of this paper reads as follows:
\medskip

\begin{theo}\label{th1.1} Then there exist two constants $\delta,C>0$ such that for any $u_0\in \dot{FB}_{1,1}^{-1}(\mathbb{R}^3)$ satisfying the condition
\begin{equation}\label{con}
\int^\infty_0||U\cd\na U||_{\dot{FB}^{-1}_{1,1}}\dd t \cdot e^{C||u_0||^2_{\dot{FB}^{-1}_{1,1}}}\leq \delta,
\end{equation}
then \eqref{NSC} admits a unique global solution
$$u\in L^\infty(0,\infty;\dot{FB}^{-1}_{1,1}(\mathbb{R}^3))\bigcap L^1(0,\infty;\dot{FB}^{1}_{1,1}(\mathbb{R}^3)).$$
\end{theo}

\begin{coro}\label{co1.1} Assume that the initial data fulfills
\begin{eqnarray}\label{initial}
\mathrm{supp} \ \hat{u}_0(\xi)\subset\mathcal{\tilde{C}}\triangleq\big\{\xi\in\R^3: \ |\xi|\geq 1\big\},
\end{eqnarray}
then there exist a sufficiently small positive constant $\delta$ and a universal constant $C$ such that if
\begin{align}\label{con1}
||u_0||_{\dot{FB}^{-1}_{1,1}}\big(||u^1_0+u^2_0,u^3_0||_{\dot{FB}^{1}_{\frac32,1}}+||\pa_3u_0||_{\dot{FB}^{1}_{\frac32,1}}\big)\cdot e^{C||u_0||^2_{\dot{FB}^{-1}_{1,1}}}\leq \delta,
\end{align}
then the system \eqref{NSC} has a unique global solution.
\end{coro}

\begin{rema}
Let two functions $a(x_1,x_2)$ with $\hat{a}(x_1,x_2)=\hat{a}(-x_1,-x_2)$ and $b(x_3)$ with $\hat{b}(x_3)=\hat{b}(-x_3)$ satisfying $\hat{a},\hat{b}\in[0,1]$,
\bbal
\mathrm{supp} \ \hat{b}\subset \{\xi_3\in \R|\ \frac12\ep<|\xi_3|<\ep\},
\end{align*}
\bbal
\hat{b}=1\quad \mathrm{on} \quad  \{\xi_3\in \R|\ \frac58\ep<|\xi_3|<\frac78\ep\},
\end{align*}
\bbal
\mathrm{supp}\ \hat{a}\subset \{\xi\in \R^2|\ |\xi_1-\xi_2|\leq \ep,\  \frac{11}{8}\leq |\xi|\leq \frac{35}{24}\},
\end{align*}
and
\bbal
\hat{a}(\xi)=1 \quad  \mathrm{on}  \quad  \{\xi\in \R^2|\ |\xi_1-\xi_2|\leq \frac12\ep,\  \frac{67}{48}\leq |\xi|\leq \frac{69}{48}\}.
\end{align*}
Then, we have for all $p\in[1,\infty]$,
\bbal
||\hat{a}_{0}||_{L^p} \approx \ep^{\frac1p}, \quad ||\hat{b}_{0}||_{L^p} \approx \ep^{\frac1p}.
\end{align*}
Let us consider the initial data $u_0=(u_{0}^1,u_{0}^2,0)$ with
\bbal
u^1_{0}=\ep^{-2}(\log\log\frac1\ep)^{\frac12}\pa_2a(x_1,x_2)b(x_3),\quad u_{0}^2=-\ep^{-2}(\log\log\frac1\ep)^{\frac12}\pa_1a(x_1,x_2)b(x_3).
\end{align*}
Direct calculation shows that
\bbal
&||u_0||_{\dot{FB}^{1}_{1,1}}\approx ||u_0||_{\dot{FB}^{-1}_{1,1}}\approx (\log\log\frac1\ep)^{\frac12},
\\&||u^1_0+u^2_0||_{\dot{FB}^{1}_{\frac32,1}}+||\pa_3u_0||_{\dot{FB}^{1}_{\frac32,1}}+||u^3_0||_{\dot{FB}^{1}_{\frac32,1}}\approx \ep^{\frac13}(\log\log\frac1\ep)^{\frac12}.
\end{align*}
Then, we can show that the left side of \eqref{con1} becomes
\begin{align*}
C\ep^{\frac13}\big(\log\log \frac1\ep\big)\exp\big(C\log\log \frac1\ep\big),
\end{align*}
which implies \eqref{NSC} have a global solution  for $\ep$ sufficiently small. For small enough $\ep$, we can deduce that $supp \ \hat{u}_0\in \{\xi\in\R^3|\ \frac43\leq |\xi|\leq \frac32\}$ and
\bbal
\De_ju_0=0,j\neq 0; \quad \De_0u_0=u_0.
\end{align*}
Therefore, we can conclude that for any $s\in\R$ and $1\leq p,r\leq \infty$
\bbal
&||u_0||_{\dot{FB}^{s}_{p,r}}\gtrsim ||u_0||_{\dot{FB}^{-1}_{1,\infty}}\approx||\hat{u}_0||_{L^1}\approx \log\log \frac1\ep.
\end{align*}
\end{rema}

\section{Proof of the main results}

{\bf Proof of Theorem \ref{th1.1}}\quad Introduce the quantity $u=U+v$, we can show that $v$ satisfies the following Cauchy problem:
\bbal\bca
\pa_tv-\De v+\Omega e_3\times v+v\cd \na v+\na p''=-U\cd \na U-v\cd\na U-U\cd \na v,
\\ \D v=0, \qquad v(0,x)=0.
\eca\end{align*}
By the Duhamel principle, this problem is equivalent to
the integral equation
$$v(t)=-\int^t_0T_\Omega(t-\tau)\mathbb{P}\big[U\cd \na U+v\cd\na v-v\cd\na U-U\cd \na v\big]d\tau,$$
where $\mathbb{P}=(\delta_{ij}+R_iR_j)_{1\leq i,j\leq3}$ denotes the Helmholtz projection onto the divergence free
vector fields, and $\{T_\Omega(t)\}_{t\geq 0}$ denotes the Stokes-Coriolis semigroup given explicitly in \cite{Hieber}.

By the similar argument of Lemma 2.2 in \cite{Sun 2019}, we have for all $t\in[0,T]$ that
\bbal
&\ \ \ \ \ ||v||_{L^\infty_t(\dot{FB}^{-1}_{1,1})}+||v||_{L^1_t(\dot{FB}^{1}_{1,1})}\\\quad
&\lesssim \int^t_0||U\cd\na U||_{\dot{FB}^{-1}_{1,1}}+||v\cd\na v||_{\dot{FB}^{-1}_{1,1}}+||U\cd\na v||_{\dot{FB}^{-1}_{1,1}}+||v\cd\na U||_{\dot{FB}^{-1}_{1,1}}\dd \tau\\\quad
&\lesssim\int^t_0||U\cd\na U||_{\dot{FB}^{-1}_{1,1}}\dd \tau +||v||_{L^2_t(\dot{FB}^{0}_{1,1})}||v||_{L^2_t(\dot{FB}^{0}_{1,1})}
+\int^t_0||v||_{\dot{FB}^{0}_{1,1}}||U||_{\dot{FB}^{0}_{1,1}}\dd \tau,
\end{align*}
where we have used Remark \ref{re1.2} and the fact $||\widehat{ab}||_{L^1}\leq ||\hat{a}||_{L^1}||\hat{b}||_{L^1}$ in the last inequality. Now, we define
\bbal
\Gamma\triangleq\sup\big\{t\in[0,T^*):||v||_{L^\infty_t(\dot{FB}^{-1}_{1,1})}
+||v||_{L^1_t(\dot{FB}^{1}_{1,1})}\leq \eta\ll1\big\},
\end{align*}
where $\eta$ is a small enough positive constant which will be determined later on. Then, it yields
\bbal
||v||_{L^\infty_t(\dot{FB}^{-1}_{1,1})}+||v||_{L^1_t(\dot{FB}^{1}_{1,1})}&\leq C\int^t_0||U\cd\na U||_{\dot{FB}^{-1}_{1,1}}\dd \tau+C\int^t_0||v||_{\dot{FB}^{-1}_{1,1}}||U||^2_{\dot{FB}^{0}_{1,1}}\dd \tau.
\end{align*}
From Gronwall's inequality, we have
\bbal
||v||_{L^\infty_t(\dot{FB}^{-1}_{1,1})}+||v||_{L^1_t(\dot{FB}^{1}_{1,1})}&\leq C\int^t_0||U\cd\na U||_{\dot{FB}^{-1}_{1,1}}\dd\tau \cdot e^{C\int^t_0||U||^2_{\dot{FB}^{0}_{1,1}}\dd \tau}
\\&\leq C\delta.
\end{align*}
Choosing $\eta=2C\delta$, thus we can get
\bbal
||v||_{L^\infty_t(\dot{FB}^{-1}_{1,1})}+||v||_{L^1_t(\dot{FB}^{1}_{1,1})}&\leq \fr\eta2 \quad\mbox{for}\quad t\leq \Gamma.
\end{align*}

So if $\Gamma<T^*$, due to the continuity of the solutions, we can obtain that there exists $0<\epsilon\ll1$ such that
\bbal
||v||_{L^\infty_t(\dot{FB}^{-1}_{1,1})}+||v||_{L^1_t(\dot{FB}^{1}_{1,1})}&\leq \eta \quad\mbox{for}\quad t\leq \Gamma+\epsilon<T^*,
\end{align*}
which is contradiction with the definition of $\Gamma$.

Thus, we can conclude $\Gamma=T^*$ and
\bbal
||v||_{L^\infty_t(\dot{FB}^{-1}_{1,1})}&\leq C<\infty \quad\mbox{for all}\quad t\in(0,T^*),
\end{align*}
which implies that $T^*=+\infty$.

%\int^t_0||U\cd\na U||_{\dot{FB}^{-1}_{1,1}}\dd\tau \cdot e^{C||u_0||^2_{\dot{FB}^{-1}_{1,1}}}
%\\&\leq C\int^t_0||U\cd\na U||_{\dot{FB}^{0}_{\frac32,1}}\dd\tau \cdot e^{C||u_0||^2_{\dot{FB}^{-1}_{1,1}}}

{\bf Proof of Corollary \ref{co1.1}}\quad Notice that $\D U$=0, we have
\bbal
&U\cd\na U^1=(U^1+U^2)\pa_1U^1+U^2\pa_2(U^1+U^2)+U^2\pa_3U^3+U^3\pa_3U^1,
\\&U\cd\na U^2=(U^1+U^2)\pa_2U^2+U^1\pa_1(U^1+U^2)+U^1\pa_3U^3+U^3\pa_3U^2,
\\&U\cd\na U^3=U^1\pa_1U^3+U^2\pa_2U^3-U^3(\pa_1U^1+\pa_2U^2).
\end{align*}
Using the fact $||ab||_{\dot{FB}^{0}_{\frac32,1}}\leq ||a||_{\dot{FB}^{0}_{\frac32,1}}||b||_{\dot{FB}^{0}_{1,1}}$, we have
\bbal
\int^t_0||U\cd\na U||_{\dot{FB}^{-1}_{1,1}}\dd\tau
&\lesssim \int^t_0||U\cd\na U||_{\dot{FB}^{0}_{\frac32,1}}\dd \tau
\\&\lesssim \int^t_0||U^1+U^2,U^3||_{\dot{FB}^{0}_{\frac32,1}\cap \dot{FB}^{1}_{\frac32,1}}||U^1,U^2||_{\dot{FB}^{0}_{1,1}\cap \dot{FB}^{1}_{1,1}}\dd\tau
\end{align*}
From \eqref{1.2}, the direct calculation shows that
\bbal
&|\hat{U}^1(\xi)|+|\hat{U}^2(\xi)|\leq e^{-t|\xi|^2}|\hat{u}_0(\xi)|,
\\&|\hat{U}^3(\xi)|\leq \Omega te^{-t|\xi|^2}\frac{|\xi_3|}{|\xi|}|\hat{u}^h_0(\xi)|+e^{-t|\xi|^2}|\hat{u}^3_0(\xi)|,
\\&|\hat{U}^1(\xi)+\hat{U}^2(\xi)|\leq te^{-t|\xi|^2}|\hat{u}^1_{0}+\hat{u}_{0}^2|+\Omega te^{-t|\xi|^2}\frac{|\xi_3|}{|\xi|}|\hat{u}_0|.
\end{align*}
This along with the property \eqref{initial} yield
\bal\label{de-con}
&\qquad \int^\infty_0||U^1+U^2,U^3||_{\dot{FB}^{0}_{\frac32,1}\cap \dot{FB}^{1}_{\frac32,1}}||U^1,U^2||_{\dot{FB}^{0}_{1,1}\cap \dot{FB}^{1}_{1,1}}\dd t
\\&\lesssim ||u_0||_{\dot{FB}^{-1}_{1,1}}\big(||u^1_0+u^2_0||_{\dot{FB}^{1}_{\frac32,1}}+||\pa_3u_0||_{\dot{FB}^{1}_{\frac32,1}}+||u^3_0||_{\dot{FB}^{1}_{\frac32,1}}\big).\nonumber
\end{align}
Thus, \eqref{de-con} is ensured whenever \eqref{con1} holds. We complete the proof of Corollary \ref{co1.1}.

\vspace*{1em}
\noindent\textbf{Acknowledgements.}  J. Li is supported by the National Natural Science Foundation of China (Grant No.11801090). J. Sun's work is partially supported by the National Natural Science Foundation of China (Grant No. 11571381, 11601434), the Natural Science Foundation of Gansu Province for Young Scholars (Grant No. 18JR3RA102).
M. Yang's work is partial supported by the National Natural Science Foundation of China (Grant No. 11801236).
%\vspace*{1em}

\end{document}